**Jointly Recovering Destination Attractiveness and Spatial Preference Segments from Observed Tours via Inverse Optimization**

Hiroyuki Hasada [a*], Yudai Honma [a]

[a] Institute of Industrial Science, The University of Tokyo, Komaba 4-6-1, Meguro-ku, Tokyo 153-8505, Japan

[*] Corresponding Author: Hiroyuki Hasada
Institute of Industrial Science, The University of Tokyo
Komaba 4-6-1, Meguro-ku, Tokyo 153-8505, Japan
Email: hasada@iis.u-tokyo.ac.jp
ORCID iD: 0000-0001-7200-7322

## AUTHOR CONTRIBUTIONS

The authors confirm contribution to the paper as follows: Conceptualization, Data curation, Methodology, Software, Visualization, Writing – original draft, Writing – review & editing: Hiroyuki Hasada; Conceptualization, Supervision, Writing – review & editing: Yudai Honma. All authors read and approved the final manuscript.

## STATEMENTS AND DECLARATIONS

### Ethical considerations
Not applicable

### Consent to participate
Not applicable

### Consent for publication
Not applicable

### Declaration of conflicting interest
The authors declare that they have no known competing financial interests or personal relationships that could have appeared to influence the work reported in this paper.

### Funding statement
This work was supported by JSPS KAKENHI [grant number JP26K17485].

### Data availability
The data used in this study were obtained from the 2021 Road Traffic Census origin–destination (OD) survey, provided by the Road Bureau of the Ministry of Land, Infrastructure, Transport and Tourism of Japan under a data-use agreement. Due to restrictions under this agreement, the data cannot be made publicly available by the authors.

## AUTHOR BIOGRAPHIES
Hiroyuki Hasada is a Research Associate at the Institute of Industrial Science, The University of Tokyo. He received his Ph.D. in Engineering from The University of Tokyo in 2024. His research develops operations research approaches using mobility-data analysis to address practical problems in transportation, urban, and architectural planning. His topics range from destination attractiveness, on-demand mobility services, and route outlier detection to inverse optimization, facility location problems, and discrete choice modeling. He also builds open-source tools and datasets for reproducible urban mobility research.

Yudai Honma is an Associate Professor at the Institute of Industrial Science, The University of Tokyo, where he leads a research laboratory. He received his Ph.D. in Engineering from Keio University in 2007. Within the institute, he also directs the Center for Social Complex Systems, which focuses on financial markets, and serves as Vice-Director of the Harmonic Mobility Research Center (ITS Center) for intelligent transport systems. He investigates spatial structures and flows across a broad range of complex systems using operations research.

# Jointly Recovering Destination Attractiveness and Spatial Preference Segments from Observed Tours via Inverse Optimization

## Abstract

Quantifying destination attractiveness by counting visitors, the most direct approach, can favor accessible sites and understate destinations that travelers reach at considerable cost. Preferences may also cut across a priori socio-demographic groupings. This study develops an inverse optimization framework that recovers destination attractiveness and preference segments jointly from observed multi-stop trip tours and generalized travel costs, without prespecifying tour-level segment assignments or imposing a parametric distribution on the unexplained part of utility. Representing each destination visit as a link in an extended network turns utility-maximizing tour choice into an inverse shortest paths problem. We apply the framework to 2,931 holiday vehicle trip tours over 51 target municipalities in Hokkaido, Japan. Constraint generation and alternating optimization make estimation tractable in this multi-destination setting, and the segments recovered under a three-segment specification show contrasting spatial preferences across the established coastal centers, the corridor around Sapporo, and the inland lake districts. The recovered values are expressed on a common monetary scale across destinations and segments without segment-specific scale normalization. They are uncorrelated with visit counts in the input data ($r = -0.01$) yet moderately correlated with tourist arrivals from an independent prefectural survey ($r = 0.40$), and are positively associated with the cost of reaching a destination ($r = 0.60$); destinations above this cost–value relationship include niche places that visit counts alone do not distinguish. The results provide a spatially explicit account of destination value and preference heterogeneity that can inform urban and regional planning and segment-targeted marketing.

## 1. INTRODUCTION

### 1.1 Background and Motivation

Quantifying the attractiveness of a place often begins with counting the number of visitors, as in destination statistics used for spatial planning. However, this straightforward approach tends to undervalue niche destinations that are difficult to access, such as remote sites that are highly valued by a small number of visitors. This distinction matters when comparing destinations with varying levels of accessibility. An alternative is therefore to evaluate destinations based on travel costs incurred by visitors (Freeman et al., 2014; Haab and McConnell, 2002). The zonal travel cost method, for example, partitions space into zones and relates zonal visit rates to travel costs to derive a demand curve, thereby valuing a site through the cost visitors are willing to bear. However, these methods typically assume a single destination and do not represent the behavior of sequentially visiting multiple destinations from a wide range of possible choices (Parsons, 2017). These approaches are also limited in representing the decision process underlying travel behavior, in which visits arise from choices made on the perceived attractiveness of destinations. This implies that attractiveness can be inversely estimated from observed travel behavior toward preferred destinations.

In decision-making for destinations, preferences are not homogeneous across individuals, and the groupings that matter for spatial choice need not be those marked by observable attributes. Socio-demographic variables describe *who holds a preference* rather than define *what that preference is*. When segments are fixed a priori by age or occupation, the analyst imposes the partition assumed to be relevant and cannot recover groups that cut across those categories, such as travelers of different ages and occupations who share the same spatial preference. Segments are better recovered from observed behavior first and characterized by socio-demographic attributes afterwards, so that the two are used in sequence rather than as substitutes.

### 1.2 Contributions

Taking observed trip tours and the travel costs between destinations as its principal empirical inputs, the proposed method jointly recovers segment-specific destination attractiveness and preference segments, without being tied to a particular spatial scale. The methodology addresses an urban and regional planning problem of comparing destinations despite differences in accessibility and visitation, while allowing relevant preference groups to emerge from behavior. As an inverse optimization approach, it empirically determines model parameters that rationalize observed decision-making behaviors (Chan et al., 2025; Chow and Recker, 2012). Here the approach takes the form of an inverse shortest paths problem (ISPP), which assigns a cost to each link such that the observed paths are the shortest in the network. We construct an extended network in which the visit to a destination is itself a link, so that the attractiveness of that destination is the cost assigned to that link, and each observed trip tour becomes a path leading from home through the visited destinations and back. Recovering these destination costs yields attractiveness values that rationalize observed path choices, consistent with the "principle of least effort" (Zipf, 1949). The recovered values thus provide a spatially explicit measure of destination attractiveness.

The methodological contribution of the proposed formulation lies in two decision variables that absorb the diversity of preferences: an error term relaxing the shortest-path constraints, and a variable for tour-level segment assignments. No single set of link costs can render every observed tour shortest when preferences differ. The error term therefore provides slack in the utility-maximization condition, with the required relaxation measuring the tour optimality gap, namely how far an observed tour falls short of optimality under the attractiveness values recovered for its segment. As a decision variable rather than a draw from an assumed distribution, it imposes no distributional form. The segment variable, in turn, groups observed tours into preference segments based on revealed behavior rather than prespecified traveler attributes. Such a non-parametric treatment of preference variation remains rare (Fosgerau, 2006).

We demonstrate the applicability of the proposed method using vehicle trip tours in Hokkaido, Japan, a region in which destinations differ widely in accessibility. To address the combinatorial number of alternative tours and the coupling of attractiveness and segment assignments, we combine constraint generation with alternating optimization, enabling region-wide estimation. The resulting maps reveal spatial preference differences among behaviorally recovered segments, distinguishing frequently visited destinations from niche destinations whose monetary attractiveness is revealed by relatively few visits despite higher travel costs. Such evidence informs how public agencies allocate tourism and recreation resources among destinations and how marketing can be targeted by segment.

## 2. LITERATURE REVIEW

### 2.1 Spatial Destination Analysis

Destination choice is commonly cast as utility maximization, with early work describing a process from awareness to an evoked set (Um and Crompton, 1990). Utility-based discrete-choice models incorporate destination and traveler attributes and highlight consumer heterogeneity (Seddighi and Theocharous, 2002). Subsequent studies examine agent-based interactions (Alvarez and Brida, 2019), spatial heterogeneity in activity-destination preferences (Chen et al., 2024), and the social context of leisure destination choice (Gramsch-Calvo and Axhausen, 2025).

Destination attractiveness has been quantified from both demand and supply perspectives. Foundational frameworks describe what visitors value through ideographic, organizational, and cognitive perspectives, using perceived attributes (Lew, 1987). For example, parametric approaches have linked destination attributes to repeat visitation using logit models (Lee et al., 2009) and have used factor analysis to identify latent dimensions of attractiveness (Boivin and Tanguay, 2019; Formica and Uysal, 2006). On the supply side, destination performance has been evaluated using non-parametric frontier methods such as data envelopment analysis to assess regional efficiency (Cracolici et al., 2008; Niavis, 2020). Non-parametric machine-learning models have examined site preferences (Chang et al., 2024).

### 2.2 Inverse Optimization of Travel Behavior

A complementary line of research uses inverse optimization to recover parameters that rationalize observed travel decisions. The household activity pattern problem formulates activity and travel scheduling as a pickup-and-delivery problem with time windows, and procedures were developed to estimate its parameters from observed schedules (Recker, 1995; Recker et al., 2008). Chow and Recker (2012) formalized this as inverse optimization, recovering the objective weights and goal arrival times under which observed schedules are optimal. The approach was later applied to freight assignment, estimating link capacities from observed flows (Chow et al., 2014), and to vehicle routing, recovering the objective weights that render observed tours optimal (You et al., 2016). This line estimates parameters for calibration units whose grouping is given rather than inferred; in You et al. (2016), which incorporates node-visit utility, each firm is treated as a given unit and calibrated separately.

Closer to choice modeling, preference heterogeneity is conventionally represented by a probability distribution, whether over latent classes (Greene and Hensher, 2003) or over the coefficients of a mixed logit (McFadden and Train, 2000). A mixing family must nonetheless be selected, and the recovered distribution is sensitive to that choice, particularly when preferences are multimodal (Fosgerau, 2006; Hess et al., 2005). Inverse optimization has therefore been used to recover taste heterogeneity without imposing such a parametric form. Hong et al. (2017) recover a joint distribution over discretized taste parameters under which observed metro routes are user-optimal. Ren and Chow (2022) estimate one coefficient vector per individual within a random-utility-consistent activity-scheduling framework. In both, the units carrying preference parameters are specified before estimation. Ren et al. (2025) use market-level share data with markets defined ex ante by population segment and origin–destination pair, and further incorporate iterative k-means clustering of market-specific attribute-coefficient vectors. Individual-to-segment assignment is not formulated as a decision variable, and the recovered quantities are attribute coefficients rather than destination-specific monetary attractiveness.

The present study lies at the intersection of these two lines, using inverse optimization to rationalize constrained multi-stop tour choices while recovering preference heterogeneity without a parametric mixing distribution. The key distinction is that parameter sharing is endogenous: tour-level segment assignments and the destination-attractiveness vector shared within each segment are estimated jointly. Because generalized travel cost enters utility with a fixed coefficient of $-1$, the recovered destination attractiveness is anchored to a common monetary scale across segments without requiring segment-specific scale normalization, whereas utility coefficients in random utility models are identified only relative to the scale of the error term. Heterogeneity is thus expressed directly in place-specific attractiveness, while residual variation is represented by endogenous tour optimality gaps rather than draws from an assumed error distribution (see Supplementary Section S6 for an empirical comparison of gap distributions). Representing destination visits as links transforms multi-stop tour choice into path choice and operationalizes this intersection as an ISPP for joint recovery of destination values and latent groups.

### 2.3 ISPP (Inverse Shortest Paths Problem)

ISPP reverses the shortest paths problem (SPP): whereas the SPP identifies the minimum-cost path by comparing alternatives, the ISPP assigns costs to each link such that every observed path is of minimum cost among its alternatives. Most ISPP formulations minimize the deviation of the estimated link costs from given initial costs, in absolute (Chen et al., 2021; Zhang et al., 1995; Zhou et al., 2014) or squared (Burton and Toint, 1992) form, subject to constraints requiring each given path to be shortest. Since the alternatives between each source and target are too numerous to enumerate explicitly, Zhang et al. (1995) generate the necessary constraints with shortest path searches while keeping the remainder implicit. Generating these constraints corresponds to column generation in the dual and constraint generation in the primal formulation used here.

Requiring every given path to be shortest is often infeasible (Bley, 2007; Call and Holmberg, 2011; Chan et al., 2025), particularly for observed travel behavior, where preferences differ across travelers. The constraints are therefore relaxed by minimizing the sum of the deviations between the given and the shortest path costs (Hasada et al., 2026; Moghaddass and Terekhov, 2021) or their maximum (Faragó et al., 2003). Heterogeneity has also been accommodated within the cost structure itself: Xu et al. (2018) assign separate link costs to each path and penalize their deviation from the average, whereas Zhou et al. (2014) treat link costs as uncertain under an assumed probability distribution.

We cast the ISPP as the inverse problem of a utility-maximizing tour choice. The deviation that existing studies minimize under the relaxed constraints is reinterpreted as an error term in utility, whose required magnitude represents the relaxation needed to satisfy tour optimality. A key extension of the ISPP is to introduce tour-level segment assignments as decision variables, enabling the joint recovery of preference segments and their shared destination attractiveness.

## 3. METHODOLOGY

### 3.1 Tour Utility Model

This section defines the tour utility model whose inverse is solved in the remainder of Section 3. The utility of a trip tour is represented as the sum of the attractiveness of the visited destinations minus the travel cost spent to visit them, together with an error term. Let $n$ $(\in N)$ index the tour records, $J_n$ denote the set of tour alternatives for record $n$, and $j$ $(\in J_n)$ index a tour alternative. Let $K$ denote the set of preference segments, with $|K|$ specified in advance, and $k(n)$ $(\in K)$ the segment assigned to record $n$, which is recovered rather than given. Let $A$ denote the destination set, $A_{nj}$ $(\subset A)$ the set of destinations visited on alternative $j$, $x_a$ the attractiveness of destination $a$ $(\in A)$ , $C_{nj}$ the generalized travel cost of alternative $j$, and $\varepsilon_{nj}$ the utility error term, which is treated not as a random variable drawn from an assumed distribution but as a decision variable recovered from the observed tours. The utility $U_{nk(n)j}$ is expressed in **Equation (1)** by separating $x'_{k(n)a}$, the portion of the attractiveness that can be distinguished by segments from the error term $\varepsilon_{nj}$, which captures residual variation not represented by the segment-specific attractiveness. By adding $x_a$ and $x'_{k(n)a}$, the attractiveness of destination $a$ for segment $k(n)$ is defined as $x_{k(n)a}$.

$$
\begin{aligned}
U_{nk(n)j} &= \textstyle\sum_{a\in A_{nj}} x_a - C_{nj} + \varepsilon_{nj}, \\
&= \textstyle\sum_{a\in A_{nj}} x_a - C_{nj} + \sum_{a\in A_{nj}} x'_{k(n)a} + \varepsilon'_{nk(n)j}, \\
&= \textstyle\sum_{a\in A_{nj}} x_{k(n)a} - C_{nj} + \varepsilon'_{nk(n)j}.
\end{aligned} \tag{1}
$$

Tour choice is assumed to maximize utility. For each record $n$ assigned to preference segment $k(n)$, the observed tour $j(n)$ is assumed to offer the greatest utility among the available alternatives in $J_n$:

$$U_{nk(n)j(n)} = \max\left\{U_{nk(n)j'} \middle| \forall j' \in J_n\right\}, \tag{2a}$$

$$\Leftrightarrow \quad U_{nk(n)j(n)} \geq U_{nk(n)j'} \qquad \forall j' \in J_n, \tag{2b}$$

$$\Leftrightarrow \quad \textstyle\sum_{a\in A_{nj(n)}} x_{k(n)a} - C_{nj(n)} + \varepsilon'_{nk(n)j(n)} \geq \sum_{a\in A_{nj'}} x_{k(n)a} - C_{nj'} + \varepsilon'_{nk(n)j'} \; \forall j' \in J_n . \tag{2c}$$

The error term of an alternative serves as a relaxation term in this maximization condition, with the error term of the revealed choice fixed at zero:

$$\varepsilon'_{nkj(n)} = 0 \qquad \forall n \in N, \forall k \in K. \tag{3}$$

When an alternative has higher utility excluding the error term than the revealed tour, a negative value of $\varepsilon'_{nk(n)j'}$ is required; its magnitude represents the relaxation required for the revealed tour to satisfy the utility-maximization condition against that alternative. Given the revealed choice result $j(n)$, the proposed model jointly recovers segment-specific destination attractiveness and segment membership by minimizing the sum of the tour optimality gaps under the assigned segments.

### 3.2 Transformation to ISPP

The above model can be transformed into an ISPP by associating destination attractiveness and travel costs with link costs on extended networks where tours are conducted. These networks are developed as follows.

First, a base network is constructed by connecting any destination or homebase with bidirectional links. Second, each site is expanded by replacing the node representing it with a link connecting its arrival and departure nodes, hereafter called an extended link. A tour is thereby represented as a directed path from the departure node of the homebase to its arrival node, in which the extended link of destination $a$ is traversed once if and only if $a$ is visited, so that its attractiveness enters the path cost once per visit, irrespective of the position at which the destination appears in the sequence. Only the travel links depend on the sequence, and their costs are given as inputs, whereas the attractiveness variables are invariant to the order of the visits. This representation corresponds to a profitable tour problem, in which a traveler selects which vertices to visit so as to minimize travel costs minus collected profits (Feillet et al., 2005). Third, the cost of each link $y_{nkl}$ is related to the attractiveness of segment $k$ or to the travel cost:

$$y_{nkl} = \begin{cases} -x_{ka} & \text{(if } l \text{ is an extended link of } a\text{)}, \\ c_{na_1a_2} & \text{(if } l \text{ is a travel link from } a_1 \text{ to } a_2\text{)}, \end{cases} \quad \forall l \in L_n, \forall k \in K, \forall n \in N, \tag{4}$$

$$x_{ka} \geq 0 \quad \forall k \in K, \forall a \in A\ . \tag{5}$$

The resulting network $G_{nk}(V, L_n, Y_{nk})$, where $V$, $L_n$ and $Y_{nk}$ denote the node set, link set, and link-cost vector, respectively, is acyclic under a distance-based ordering of zones, with destination attractiveness constrained to be non-negative (see Supplementary Section S2). The alternative set $J_n$ consists of all tours represented by directed paths from the tour's start node to its end node on this network, including the observed tour $j(n)$ but excluding visit sequences inconsistent with the zone order.

On the resulting network, the cost of path $j$ comprising links $L_{nj}$, $\sum_{l \in L_{nj}} y_{nkl}$, corresponds to the additive inverse of the utility excluding the error term of tour alternative $j$ that visits the destinations between the start and end of the tour. Consequently, as shown below, using a path with a lower cost is equivalent to selecting a tour with higher utility.

$$\sum_{l \in L_{nj(n)}} y_{nkl} = \min\left\{\sum_{l \in L_{nj'}} y_{nkl} - \varepsilon'_{nkj'} \middle| \forall j' \in J_n\right\}, \quad \forall n \in N, \forall k \in K, \tag{6a}$$

$$\Leftrightarrow \quad \sum_{l \in L_{nj(n)}} y_{nkl} \leq \sum_{l \in L_{nj'}} y_{nkl} - \varepsilon'_{nkj'} \quad \forall j' \in J_n, \forall n \in N, \forall k \in K, \tag{6b}$$

$$\Leftrightarrow \quad -\sum_{a \in A_{nj(n)}} x_{ka} + C_{nj(n)} \leq -\sum_{a \in A_{nj'}} x_{ka} + C_{nj'} - \varepsilon'_{nkj'} \forall j' \in J_n, \forall n \in N, \forall k \in K, \tag{6c}$$

$$\Leftrightarrow \quad U_{nkj(n)} = \max\left\{U_{nkj'} \middle| \forall j' \in J_n\right\}, \quad \forall n \in N, \forall k \in K, \tag{6d}$$

where **Equation (3)** holds. Thus, decision making based on the utility of visiting destinations is equivalent to selecting the shortest path on network $G_{nk}(V, L_n, Y_{nk})$. Given the observed paths, the ISPP determines each link cost in a network such that each path cost is the minimum of the alternatives. By constraining **Equations (6)** with the input of the paths corresponding to the selected tours and fixing the travel costs, the link costs distinguished by each segment can be inversely estimated, deriving the attractiveness of the destinations for each segment. The introduced error term, $\varepsilon'_{nkj}$, functions as a relaxation term for the shortest-path constraint in the ISPP. This approach relaxes the impractical constraint that all given paths be shortest under a common set of link costs. This allows for heterogeneity in preferences.

### 3.3 Formulation

The relaxed shortest path constraint (**Equation (6b)**) is imposed with respect to the link costs for each preference segment. In the proposed formulation, a binary variable $r_{nk}$ is introduced, indicating whether record $n$ is assigned to segment $k$. This results in a mixed-integer bilinear programming formulation:

$$\min_{x,r,\delta,\varepsilon'} \quad \sum_{n,k} r_{nk}\delta_{nk}. \tag{7}$$

s.t.

**Constraints (3), (4), (5), (6b)**

$$r_{nk} \in \{0,1\} \quad \forall n \in N, \forall k \in K. \tag{8}$$

$$\sum_k r_{nk} = 1 \quad \forall n \in N. \tag{9}$$

$$\varepsilon'_{nkj'} \geq -\delta_{nk} \quad \forall j' \in J_n, \forall n \in N, \forall k \in K. \tag{10}$$

**Objective function (7)** minimizes the sum of tour optimality gaps under the assigned segments. **Constraint (8)** defines $r_{nk}$ as a binary variable equal to 1 if record $n$ is assigned to segment $k$, and 0 otherwise, indicating segment membership. **Constraint (9)** ensures that each record is assigned to exactly one segment. **Constraint (10)** bounds $-\varepsilon'_{nkj'}$ by $\delta_{nk}$ over all alternative tours. Consequently, only an alternative with a lower path cost than the observed tour imposes a positive lower bound on $\delta_{nk}$; if the observed tour is already a shortest path, the gap is zero.

Two main types of decision variables are obtained: $r_{nk}$ for tour-level segment assignments and $x_{ka}$ for segment-specific destination attractiveness. The error term $\varepsilon'_{nkj}$ and the gap variable $\delta_{nk}$ are auxiliary variables used to formulate the relaxed optimality condition.

At an optimum, $\delta_{nk}$ may be taken at its smallest feasible value under **Constraints (6b) and (10)**, which equals the gap between the observed path cost and the shortest path cost under the link costs of segment $k$:

$$\delta_{nk} = \sum_{l \in L_{nj(n)}} y_{nkl} - \min\left\{\sum_{l \in L_{nj'}} y_{nkl} \middle| \forall j' \in J_n\right\}. \tag{11}$$

For the segment assigned to record $n$, **Objective function (7)** drives $\delta_{nk}$ to this value; for other segments, taking $\delta_{nk}$ at its smallest feasible value does not affect the objective. Because the resulting gap depends only on the link costs, it is well defined for every segment, including those to which record $n$ is not assigned. This property is used to examine overlapping membership in the following section.

### 3.4 Solution Algorithm

The proposed model poses two computational challenges: enumerating all possible paths generates an excessive number of constraints, and the objective function involves a mixed-integer bilinear formulation. The proposed approach combines constraint generation with alternating optimization. Because the model is an ISPP, constraint generation identifies and adds violated constraints through shortest-path searches on the extended network, making region-wide estimation tractable without enumerating alternative tours. Alternating optimization updates link costs and segment-assignment variables in turn after relaxing the latter from binary to continuous, as in k-means clustering (MacQueen, 1967). Because alternating optimization may converge to a local solution, global optimality is not guaranteed. The algorithm also identifies observed tours equally well explained by multiple segments based on their tour optimality gaps. **Figure 1** outlines the algorithm, with detailed procedures and pseudocode provided in Supplementary Section S3.

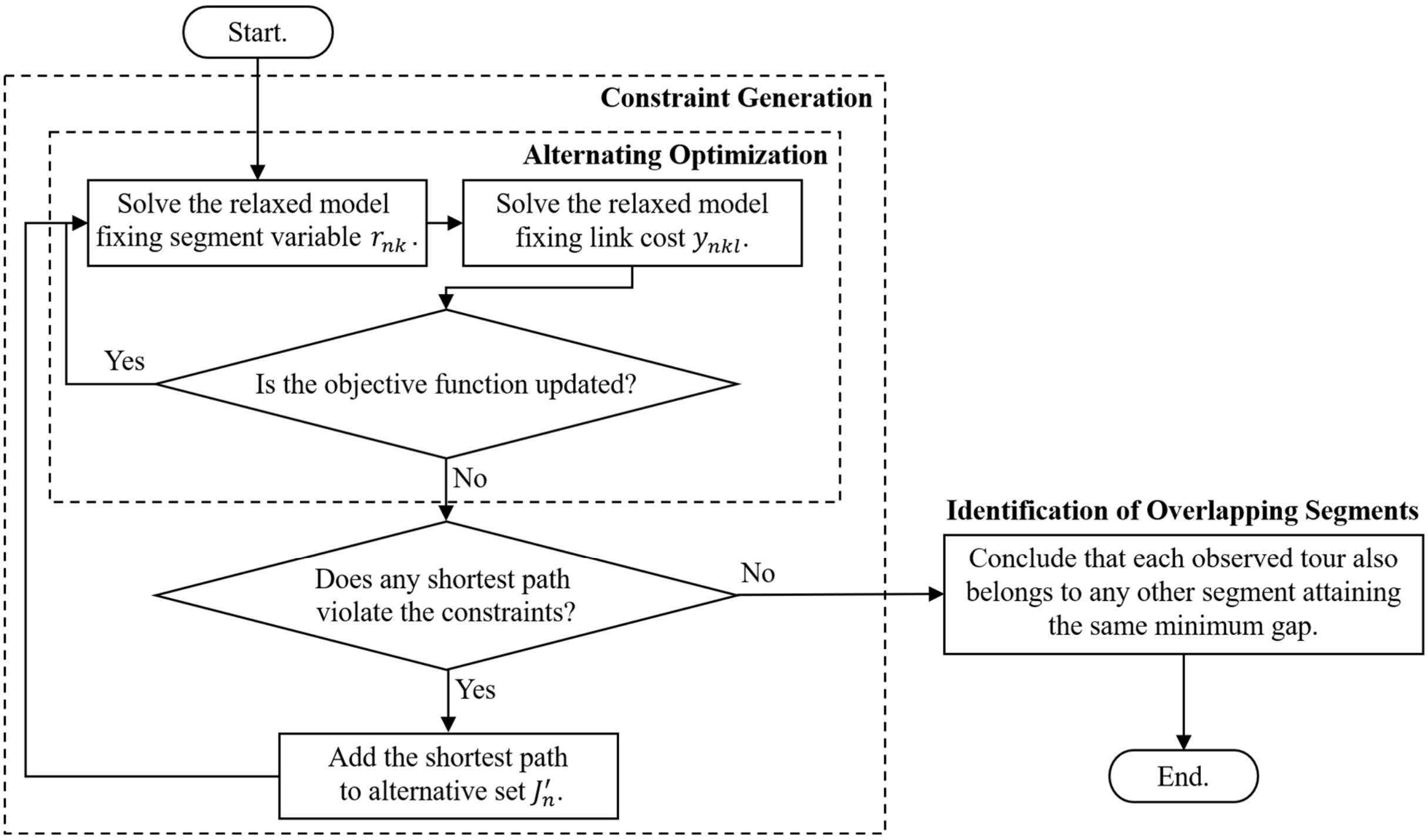

**Figure 1 Outline of the proposed solution algorithm**

## 4. STUDY AREA AND DATA

This study used Japan's national origin–destination survey, which recorded household and vehicle attributes and trip details in the autumn of 2021, aggregated at the municipal zone level. The analysis focuses on weekend and holiday vehicle trips within Hokkaido, Japan's northernmost prefecture, where leisure resources are dispersed across a sparsely populated area and access is almost entirely car-based. Consecutive leisure-related trips were linked into tours. The resulting input tours contain neither everyday activities nor sub-loops, and visit each destination at most once.

Travel link costs are specified as constants common to all tours for this empirical application, calculated as the sum of the opportunity cost of travel time and fuel cost. The opportunity cost of travel time was set to one-third of the value of travel time (Kubo et al., 2020; Sinclair et al., 2020). The extended link costs are defined as the negatives of the destination attractiveness variables to be recovered; a separate variable was assigned to each of the 51 destinations visited at least 20 times in the input data, hereafter the target destinations, while the other destinations share a single common variable to limit the solution space and improve identifiability.

## 5. RESULTS

This section presents the segments and destination attractiveness recovered from the trip tour path data. **Table 1** summarizes the results, where $|K|$ represents the number of segments predefined in the proposed model. Segments are labeled 1A through 3C in descending order of the number of tours belonging to each segment. The values for $|K| = 1$ in the preference segmentation section represent the properties of the input data. It should be noted that the total number of tours in all segments may not be equal to the 2,931 tours in the input data, as the proposed method allows for overlapping segments.

**Table 1 Results of recovered segments and attractiveness**

| | | $\lvert K \rvert = 1$ | $\lvert K \rvert = 2$ | | $\lvert K \rvert = 3$ | | |
|---|---|---|---|---|---|---|---|
| | | 1A | 2A | 2B | 3A | 3B | 3C |
| **Preference segmentation** | | | | | | | |
| Counts of tours | | 2931 | 1984 | 1312 | 1434 | 1181 | 1039 |
| Counts of tours belonging to only one segment | | 2931 | 1619 | 947 | 1049 | 768 | 578 |
| Counts of destinations in each trip tour | avg. | 1.15 | 1.15 | 1.11 | 1.13 | 1.08 | 1.15 |
| Travel cost of each trip tour [yen] | avg. | 3304.3 | 3200.6 | 2938.1 | 3151.0 | 2816.7 | 2883.9 |
| **Destination attractiveness** | | | | | | | |
| Attractiveness of target destinations [yen] | avg. | 587.5 | 651.6 | 727.0 | 704.2 | 571.8 | 726.9 |
| | std. dev. | 702.5 | 912.6 | 1740.8 | 1753.3 | 1307.3 | 1311.1 |
| | max. | 3737.8 | 3814.1 | 9236.8 | 11027.2 | 7358.4 | 7191.1 |
| | min. | 0.0 | 0.0 | 0.0 | 0.0 | 0.0 | 0.0 |
| Attractiveness of other destinations [yen] | | 69.0 | 69.0 | 149.0 | 34.9 | 313.2 | 149.0 |

The sizes of the recovered segments differ: for $|K| = 2$, 1,984 tours belong to Segment 2A and 1,312 to Segment 2B. In addition, the number of tours belonging to only one segment decreased under the three-segment specification. For example, in Segment 3C, only 578 out of 1,039 tours belong exclusively to that segment. This growing overlap indicates that, as additional segments are introduced, more observed tours are explained equally well by multiple segments. The overlap ratio therefore diagnoses possible over-segmentation, and the $|K| = 3$ specification used below is exploratory rather than an identified optimal number of segments. The average number of destinations per tour was relatively constant across segments, whereas average travel costs varied, reflecting differences in the travel patterns and the areas each segment visited, which are examined in Section 6.

The model recovers segment-specific values for the same destination, and the values within each segment are widely dispersed, as indicated by the large standard deviations and wide ranges in **Table 1**. The model thus assigns a destination-specific and segment-specific value to each municipal destination. The attractiveness of destinations with no visits in any segment was generally estimated to be 0, because positive attractiveness values would only increase the utility of unchosen alternatives.

To evaluate computational performance, we compared the proposed algorithm with a brute-force baseline that enumerates all candidate paths and solves the resulting mixed-integer bilinear program directly. The comparison used 20 sampled datasets of 100 tours and was restricted to single-destination tours to make full enumeration feasible, although the proposed framework also accommodates multi-destination tours. The baseline consistently reached the 3,600-s time limit without proving optimality, whereas the proposed algorithm terminated within 20–170 s and at most 13 constraint-generation iterations. The proposed algorithm also attained lower objective values on average than the baseline's time-limit incumbents, although global optimality was not guaranteed. Detailed experimental settings and results are provided in Supplementary Section S5.

## 6. DISCUSSION

### 6.1 Empirical Characteristics of Recovered Segments

This section discusses the results when the number of segments is set to three. **Figure 2** shows the study area together with the visits each destination received and the trips made within each of the three segments in the input trip tour data, with the most populous municipal zones labeled. The segments are then examined against socio-demographic attributes. **Figure 3** illustrates the composition ratios of tours belonging to only one segment according to the traveler's age group and occupation and the vehicle type used. The proportion in Segment 3A generally increases with age, whereas that in Segment 3C decreases. Students and homemakers were primarily classified into Segments 3B and 3C, which had lower average travel costs, as shown in **Table 1**, consistent with a preference for leisure destinations closer to home. Chartered buses accounted for a higher proportion of Segment 3C than other vehicle types did; this segment also visits the most destinations per tour (**Table 1**), a combination consistent with participation in packaged group itineraries. These associations are post hoc descriptions of the observed sample and do not provide an out-of-sample rule for segment assignment.

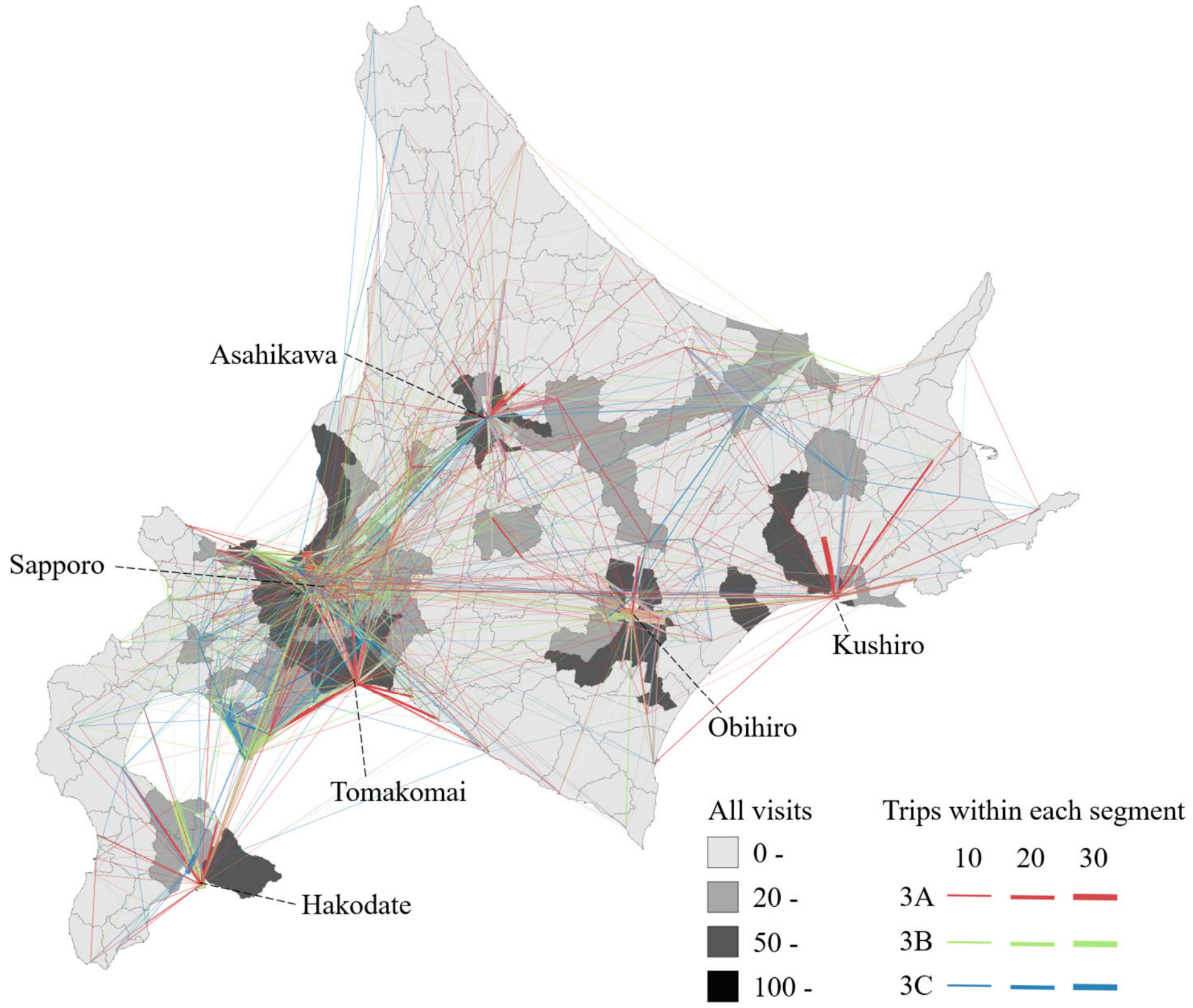

**Figure 2 Spatial distribution of visits and trips by segment among the input trip tour data**

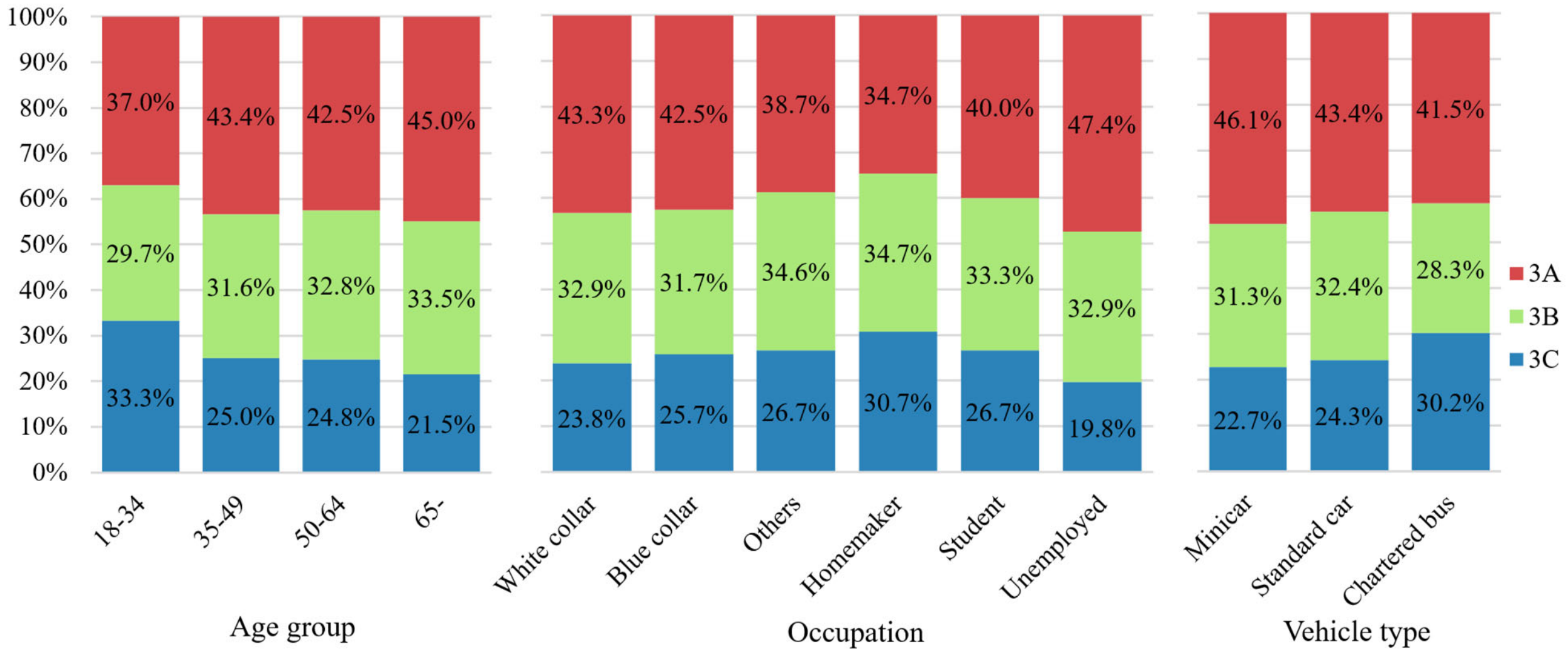

**Figure 3 Ratios of each segment by socio-demographic attributes**

## 6.2 Spatial Patterns of Recovered Attractiveness

**Figure 4** maps the highest attractiveness each destination attains across the three segments and the segment with the highest recovered value. The five highest-valued destinations in each segment are labeled unless already labeled among the most populous zones in **Figure 2**. Segment 3A values Hakodate, Kushiro and Tomakomai, among the most populous municipal zones in **Figure 2**. Segment 3B values Obihiro and Abashiri in the east, together with Furano and Otaru, which lie close to Sapporo, Hokkaido's capital city. Segment 3C combines regional centers such as Asahikawa and Kitami with lakeside destinations, such as Teshikaga, located around Lakes Kussharo and Mashu, as well as Niseko and Date, situated near Lakes Toya and Shikotsu. Read together, **Figures 2 and 4** show that the segments divide the island between the established coastal centers, the corridor around Sapporo, and the inland lake districts. Destinations within and adjacent to Sapporo attain lower values in every segment, as they are reached at little cost and need not be highly valued to rationalize the observed tours.

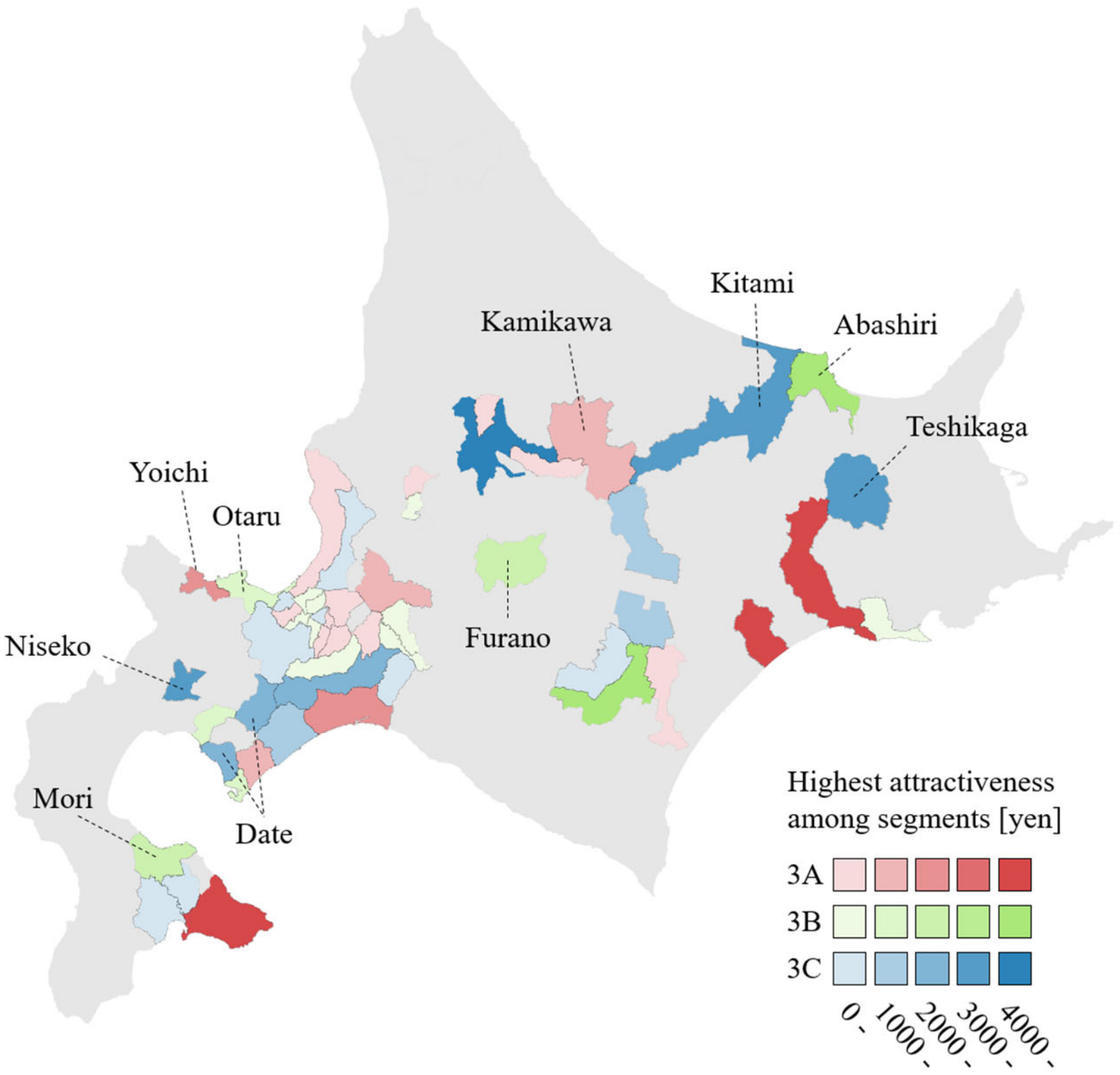

**Figure 4 Spatial distribution of the highest recovered attractiveness among segments**

The recovered values are next compared with indicators that can be observed directly (**Figure 5**). Attractiveness rises with the travel cost incurred to reach a destination from home, averaged over the tours visiting it (Pearson correlation 0.60 over the 51 target destinations), and is uncorrelated with the number of visits recorded in the input data (−0.01). It is moderately correlated with the visitor counts reported for the same season in the prefectural tourism survey (0.40 over the 41 destinations covered by the survey; Hokkaido Government, 2022). The estimates therefore do not reproduce visit frequency but reflect the travel costs incurred to reach a destination. The dashed curve in **Figure 5** is the envelope below which 90% of the destinations lie, and the destinations on or above it are those whose recovered value stands out against the cost of reaching them.

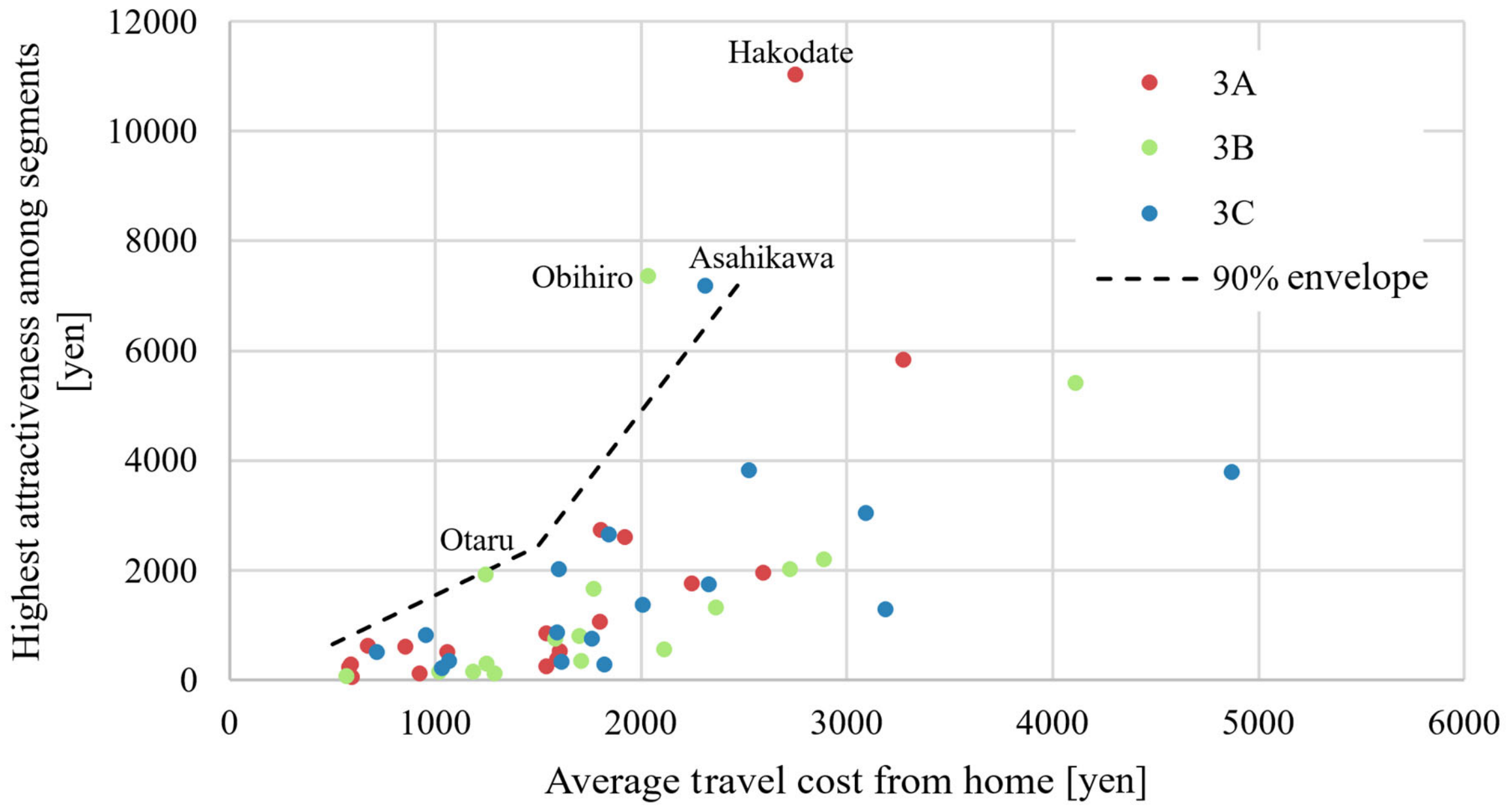

**Figure 5 Highest recovered attractiveness of each destination against the average travel cost from home**

These findings have two implications for destination planning. First, the near-zero correlation between visitation and recovered attractiveness indicates that visit counts alone may obscure destination values relative to access costs, as illustrated by Hakodate, Obihiro, Asahikawa, and Otaru in **Figure 5**. Second, **Figures 2 and 4** reveal contrasting spatial preferences across segments, from coastal centers to the corridor around Sapporo and the inland lake districts. Together, these results support spatial comparisons among destinations without assuming spatially homogeneous preferences, although they are diagnostic rather than predictive and do not establish the effects of specific planning interventions.

## 7. CONCLUSION

This study developed an inverse optimization framework that jointly recovers segment-specific destination attractiveness and preference segments from observed multi-stop tours and generalized travel costs. Applied to holiday vehicle tours in Hokkaido, the proposed ISPP identified destinations whose revealed attractiveness was high relative to access costs and mapped contrasting spatial preferences across segments. These outputs provide a spatially explicit characterization of how destination attractiveness and preference heterogeneity vary across a regional destination system. By complementing visitor counts with destination values recovered jointly from observed tours and generalized travel costs, the framework provides a quantitative screening tool for public agencies and destination-management organizations.

The utility specification intentionally focuses on destination attractiveness and generalized travel costs. Accordingly, the recovered values represent attractiveness revealed through observed tour choices; they are not intended to explain or predict attractiveness from destination attributes or other underlying mechanisms. Factors not represented separately, including such attributes, prior experience, entrance fees, and other activity-related considerations, may instead contribute to the tour optimality gaps. These boundaries should guide interpretation and use of the recovered values.

Future work should establish formal identifiability conditions and methods for assessing uncertainty in the recovered attractiveness. More systematic segment selection and validation remain future tasks beyond the present scope; observed overlap could serve as one diagnostic signal in developing such criteria. Empirical robustness should be examined across seasons, urban and rural settings, travel modes, and data sources. A complementary direction is to extend the spatial, acyclic network representation to a time-space structure, allowing activity durations, time windows, scheduling constraints, and richer feasible tour sequences to be represented explicitly and suggesting potential applications to time-space network design problems (Tong et al., 2015). Addressing the computational burden of larger state spaces and more complex tour alternatives will help broaden the framework's applicability to richer mobility data and more complex real-world settings.

**References**

Alvarez E and Brida JG (2019) An agent-based model of tourism destinations choice. *International Journal of Tourism Research* 21(2): 145–155.

Bley A (2007) Inapproximability results for the inverse shortest paths problem with integer lengths and unique shortest paths. *Networks* 50(1): 29–36.

Boivin M and Tanguay GA (2019) Analysis of the determinants of urban tourism attractiveness: The case of Québec City and Bordeaux. *Journal of Destination Marketing & Management* 11: 67–79.

Burton D and Toint PL (1992) On an instance of the inverse shortest paths problem. *Mathematical Programming* 53: 45–61.

Call M and Holmberg K (2011) Complexity of inverse shortest path routing. In: Pahl J, Reiners T and Voß S (eds) *Network Optimization*. Berlin: Springer, pp.339–353.

Chan TCY, Mahmood R and Zhu IY (2025) Inverse optimization: Theory and applications. *Operations Research* 73(2): 1046–1074.

Chang V, Islam MR, Ahad A, et al. (2024) Machine learning for predicting tourist spots' preference and analysing future tourism trends in Bangladesh. *Enterprise Information Systems* 18(12): 1450–1481.

Chen F, Zhou S, Lu J, et al. (2024) A behavioral explanation of the activity-space segregation: Individuals' preference of choosing an activity destination. *Environment and Planning B: Urban Analytics and City Science* 51(9): 1995–2011.

Chen L, Chen Y and Langevin A (2021) An inverse optimization approach for a capacitated vehicle routing problem. *European Journal of Operational Research* 295(3): 1087–1098.

Chow JYJ and Recker WW (2012) Inverse optimization with endogenous arrival time constraints to calibrate the household activity pattern problem. *Transportation Research Part B: Methodological* 46(3): 463–479.

Chow JYJ, Ritchie SG and Jeong K (2014) Nonlinear inverse optimization for parameter estimation of commodity-vehicle-decoupled freight assignment. *Transportation Research Part E: Logistics and Transportation Review* 67: 71–91.

Cracolici MF, Nijkamp P and Rietveld P (2008) Assessment of tourism competitiveness by analysing destination efficiency. *Tourism Economics* 14(2): 325–342.

Faragó A, Szentesi Á and Szviatovszki B (2003) Inverse optimization in high-speed networks. *Discrete Applied Mathematics* 129(1): 83–98.

Feillet D, Dejax P and Gendreau M (2005) Traveling salesman problems with profits. *Transportation Science* 39(2): 188–205.

Formica S and Uysal M (2006) Destination attractiveness based on supply and demand evaluations: An analytical framework. *Journal of Travel Research* 44(4): 418–430.

Fosgerau M (2006) Investigating the distribution of the value of travel time savings. *Transportation Research Part B: Methodological* 40(8): 688–707.

Freeman AM III, Herriges JA and Kling CL (2014) *The Measurement of Environmental and Resource Values: Theory and Methods*. 3rd ed. Abingdon: Routledge.

Gramsch-Calvo B and Axhausen KW (2025) The importance of the social environment on leisure destination choice: A mixed multinomial analysis of homophilic preferences. *Environment and Planning B: Urban Analytics and City Science* 52(7): 1650–1668.

Greene WH and Hensher DA (2003) A latent class model for discrete choice analysis: Contrasts with mixed logit. *Transportation Research Part B: Methodological* 37(8): 681–698.

Haab TC and McConnell KE (2002) *Valuing Environmental and Natural Resources: The Econometrics of Non-market Valuation*. Cheltenham: Edward Elgar.

Hasada H, Hasegawa D, Zhou XS, et al. (2026) Unveiling implicit path outlier deviations via foundational linear inverse optimization in mobility networks. *Transportmetrica A: Transport Science*. Epub ahead of print 16 March 2026. DOI: 10.1080/23249935.2026.2635567.

Hess S, Bierlaire M and Polak JW (2005) Estimation of value of travel-time savings using mixed logit models. *Transportation Research Part A: Policy and Practice* 39(2–3): 221–236.

Hokkaido Government (2022) Survey report on tourist arrivals in Hokkaido [in Japanese]. Available at: https://www.pref.hokkaido.lg.jp/kz/kkd/irikomi.html (accessed 15 August 2026).

Hong SP, Kim KM, Byeon G, et al. (2017) A method to directly derive taste heterogeneity of travellers' route choice in public transport from observed routes. *Transportation Research Part B: Methodological* 95: 41–52.

Kubo T, Uryu S, Yamano H, et al. (2020) Mobile phone network data reveal nationwide economic value of coastal tourism under climate change. *Tourism Management* 77: 104010.

Lee C-F, Ou W-M and Huang H-I (2009) A study of destination attractiveness through domestic visitors' perspectives: The case of Taiwan's Hot Springs tourism sector. *Asia Pacific Journal of Tourism Research* 14(1): 17–38.

Lew A (1987) A framework of tourist attraction research. *Annals of Tourism Research* 14(4): 553–575.

McFadden D and Train K (2000) Mixed MNL models for discrete response. *Journal of Applied Econometrics* 15(5): 447–470.

MacQueen J (1967) Some methods for classification and analysis of multivariate observations. In: Le Cam LM and Neyman J (eds) *Proceedings of the Fifth Berkeley Symposium on Mathematical Statistics and Probability, Vol. 1*. Berkeley, CA: University of California Press, pp.281–297.

Moghaddass M and Terekhov D (2021) Inverse integer optimization with multiple observations. *Optimization Letters* 15: 1061–1079.

Niavis S (2020) Evaluating the spatiotemporal performance of tourist destinations: The case of Mediterranean coastal regions. *Journal of Sustainable Tourism* 28(9): 1310–1331.

Parsons GR (2017) Travel cost models. In: Champ PA, Boyle KJ and Brown TC (eds) *A Primer on Nonmarket Valuation*. 2nd ed. Dordrecht: Springer, pp.187–233.

Recker WW (1995) The household activity pattern problem: General formulation and solution. *Transportation Research Part B: Methodological* 29(1): 61–77.

Recker WW, Duan J and Wang H (2008) Development of an estimation procedure for an activity-based travel demand model. *Computer-Aided Civil and Infrastructure Engineering* 23(7): 483–501.

Ren X and Chow JYJ (2022) A random-utility-consistent machine learning method to estimate agents' joint activity scheduling choice from a ubiquitous data set. *Transportation Research Part B: Methodological* 166: 396–418.

Ren X, Chow JYJ and Bansal P (2025) Nonparametric mixed logit model with market-level parameters estimated from market share data. *Transportation Research Part B: Methodological* 196: 103220.

Seddighi HR and Theocharous AL (2002) A model of tourism destination choice: A theoretical and empirical analysis. *Tourism Management* 23(5): 475–487.

Sinclair M, Mayer M, Woltering M, et al. (2020) Valuing nature-based recreation using a crowdsourced travel cost method: A comparison to onsite survey data and value transfer. *Ecosystem Services* 45: 101165.

Tong L, Zhou X and Miller HJ (2015) Transportation network design for maximizing space–time accessibility. *Transportation Research Part B: Methodological* 81: 555–576.

Um S and Crompton J (1990) Attitude determinants in tourism destination choice. *Annals of Tourism Research* 17(3): 432–448.

Xu SJ, Nourinejad M, Lai X, et al. (2018) Network learning via multiagent inverse transportation problems. *Transportation Science* 52(6): 1347–1364.

You SI, Chow JYJ and Ritchie SG (2016) Inverse vehicle routing for activity-based urban freight forecast modeling and city logistics. *Transportmetrica A: Transport Science* 12(7): 650–673.

Zhang J, Ma Z and Yang C (1995) A column generation method for inverse shortest path problems. *ZOR – Methods and Models of Operations Research* 41(3): 347–358.

Zhou J, Yang F and Wang K (2014) An inverse shortest path problem on an uncertain graph. *Journal of Networks* 9(9): 2353–2359.

Zipf GK (1949) *Human Behavior and the Principle of Least Effort*. Cambridge, MA: Addison-Wesley Press.

## SUPPLEMENTAL MATERIAL

### S1. Literature Review of Non-parametric Approaches

Data-driven non-parametric methods are widely used to analyze visit behavior without assuming a functional form, although they do not represent the choice process that generates a tour. Mixture- and tree-based predictors support prediction and segmentation (Bishop, 2006), topic models (Blei et al., 2003) and probabilistic factorization (Mnih and Salakhutdinov, 2007) uncover co-visitation structures, and frontier methods benchmark efficiency (Cooper et al., 2011). None of them encodes the combinatorial feasibility or least-cost consistency of multi-stop tours by construction or yields a structural tour-choice model. The proposed ISPP complements these approaches by instead preserving tour-level optimality over a very large alternative set.

### S2. Construction of an Acyclic Extended Network

Starting from a base network connecting all destinations and the homebase with bidirectional travel links, the extended network is constructed by splitting each site into arrival and departure nodes joined by an extended link. The cost of each destination's extended link is set to the negative of its attractiveness, which is restricted to be non-negative under the assumption that visiting a destination does not reduce tour utility. Consequently, extended links have non-positive costs, and the network may contain directed cycles of negative total cost, which would make the shortest-path cost unbounded below and allow a destination to be counted more than once. Such cycles are therefore excluded by ordering the zones according to their distance from the start of the tour and retaining only the links consistent with this order, which removes a travel link connecting a distant zone to a closer one unless its head is the end node, as well as the links originating from the end node (Figure S1). Every remaining link either advances the zone order or enters the end node, while each extended link is directed only from the arrival node to the departure node. The resulting network $G_{nk}(V, L_n, Y_{nk})$, where $V$, $L_n$ and $Y_{nk}$ denote the node set, link set, and the vector of link costs, respectively, is acyclic and its shortest path can be obtained by dynamic programming despite the negative link costs. The alternative set $J_n$ consists of all tours represented by directed paths from the start node to the end node on this network, including multi-destination tours but excluding visit sequences inconsistent with the zone order. The observed tour $j(n)$ is included in $J_n$.

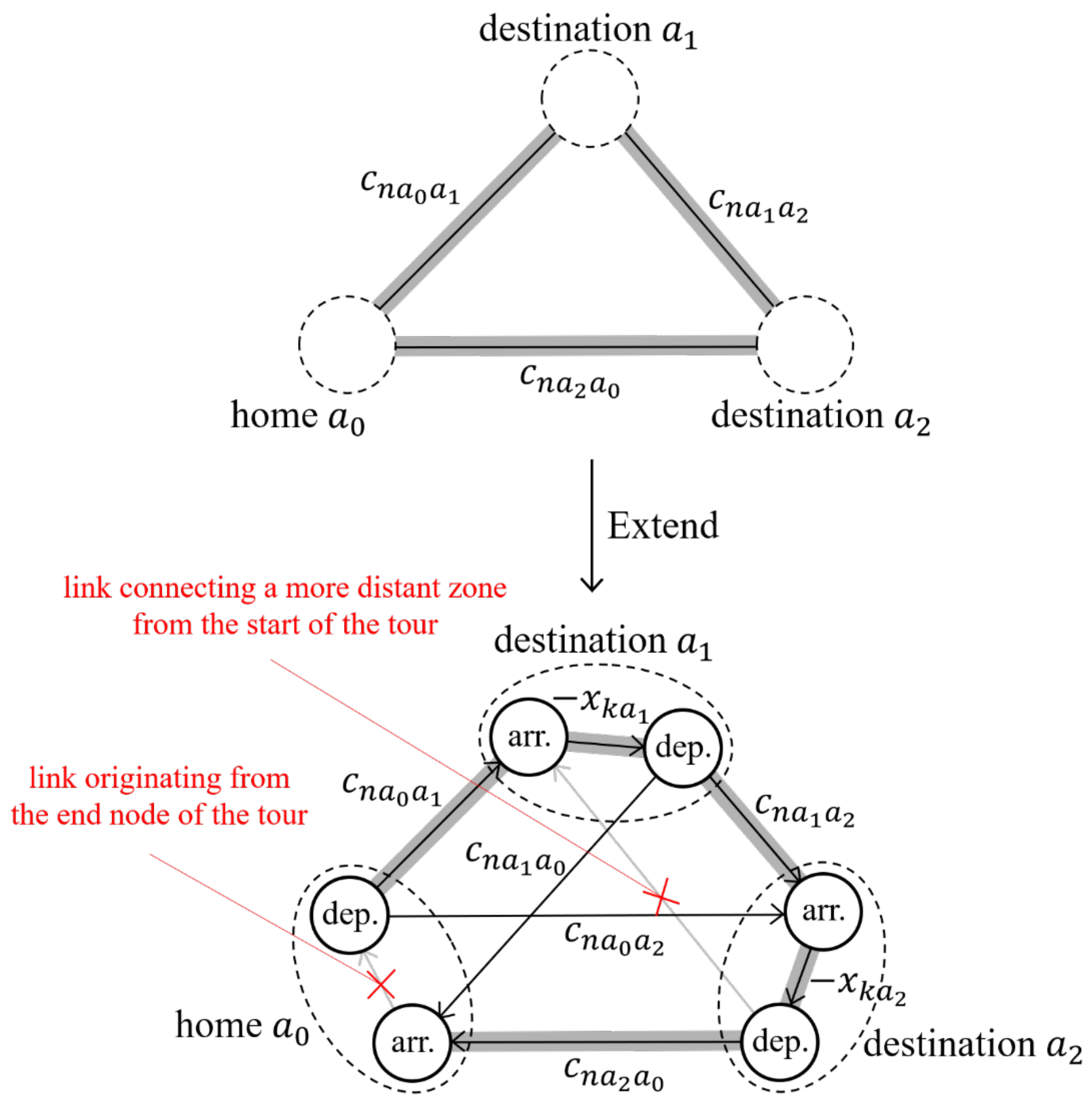

**Figure S1 Construction of the extended network and the links removed to prevent negative cycles**

## S3. Detailed Procedures and Pseudocode of the Solution Algorithm

This section provides the detailed procedures and pseudocode for the algorithm outlined in Figure 1 of the main text, using the same notation and equation numbers.

*S3.1 Constraint Generation*

Constraint generation solves a restricted model that imposes **Constraint (6b)** on a subset $J_n' \subseteq J_n$ and then adds the violated constraints, following the master–subproblem structure of column generation applied to rows. Eliminating $\varepsilon'_{nkj'}$ between **Constraints (6b) and (10)** shows that the violation of an alternative increases as its path cost decreases, so identifying the most violated constraint reduces to finding the shortest path on $G_{nk}(V, L_n, Y_{nk})$, which can be found exactly by dynamic programming because the network is acyclic. The pseudocode is presented below:

| **Algorithm S1 Constraint generation for the proposed model** |
|---|
| $J_n' \leftarrow \{j(n)\}$ |
| **while** any $J_n'$ is updated **do** |
| find the solution $y^*_{nkl}$ and $r^*_{nk}$ of the restricted model using **Algorithm S2** |
| **for each** $n \in N$ **do** |
| find the shortest path $\tilde{j}$ on $G_{nk(n)}(V, L_n, Y^*_{nk(n)})$ between start and end of tour $n$ |
| **if** min $\left\{\sum_{l\in L_{nj'}} y^*_{nk(n)l} \,\middle|\, \forall j' \in J_n'\right\} - \sum_{l\in L_{n\tilde{j}}} y^*_{nk(n)l} > \theta$ **then** |
| update $J_n' \leftarrow J_n' \cup \{\tilde{j}\}$ |
| **end if** |
| **end for** |
| **end while** |

A path already in $J'_n$ satisfies the current constraints and is never added twice, so each iteration strictly enlarges $J'_n$ and the algorithm terminates finitely on the finite set $J_n$. The restricted model is a relaxation and the shortest path search identifies the violation exactly. Given an optimal solution of each restricted model, the algorithm therefore attains a global optimum of the full model, because only the assigned segments enter **Objective function (7)**. This global-optimality result is conditional on each restricted model being solved globally; the alternating optimization used below does not provide this guarantee.

*S3.2 Alternating Optimization*

Although the number of constraints is reduced in the preceding section, solving large-scale bilinear programming problems remains challenging. The restricted model solved within the constraint generation is therefore optimized by alternating between the two blocks of variables, in the manner of the k-means clustering algorithm, a specific case of the more general expectation–maximization algorithm (Dempster et al., 1977). The proposed algorithm minimizes the sum of the product of segment variable $r_{nk}$ and the tour optimality gap $\delta_{nk}$ in **Objective function (7)**. The binary variables are first relaxed into continuous ones:

$$r_{nk} \geq 0 \qquad \forall n \in N, \forall k \in K. \tag{8'}$$

This relaxation preserves the optimal objective value: with link costs fixed, **Constraints (6b) and (10)** determine the minimum feasible gap for each record–segment pair, so each record can be assigned entirely to a minimum-gap segment. When several segments attain the same minimum gap, fractional assignments may also be optimal.

The pseudocode is presented below:

| **Algorithm S2 Alternating optimization of the restricted model** |
|---|
| $\tilde{r}_{nk}, \tilde{r}^*_{nk}$ and $\tilde{\delta}^*{}_{nk}$ are randomly set |
| $\tilde{\delta}_{nk}$ is set to a sufficiently large number |
| **while** $\sum_{n,k} \tilde{r}_{nk}\, \tilde{\delta}_{nk} - \sum_{n,k} \tilde{r}^*_{nk}\, \tilde{\delta}^*{}_{nk} > \theta$ **do** |
| **if** $\sum_n \tilde{r}^*_{nk} = 0$ for any $k$ **then** |
| $\tilde{r}^*_{nk}$ is randomly set |
| **end if** |
| update $\tilde{r}_{nk} \leftarrow \tilde{r}^*_{nk}$ and $\tilde{\delta}_{nk} \leftarrow \tilde{\delta}^*{}_{nk}$ |
| find the optimal solution $\tilde{y}^*_{nkl}$ of the linear model fixing $r_{nk}$ to $\tilde{r}_{nk}$ |
| find the optimal solution $\tilde{r}^*_{nk}$ and $\tilde{\delta}^*{}_{nk}$ of the relaxed model fixing $y_{nkl}$ to $\tilde{y}^*_{nkl}$ |
| **end while** |

Each block is optimized with the other block fixed, so **Objective function (7)** is non-increasing over the iterations and bounded below by zero, and the procedure stops when the improvement falls below the near-zero threshold $\theta$. Empty segments are reset. As in k-means, the procedure yields a local solution of the restricted model and may depend on initialization.

*S3.3 Identification of Overlapping Segments*

**Constraint (9)** assigns each record to a single segment, whereas the assignment that minimizes **Objective function (7)** need not be unique. After the constraint generation terminates, the gap $\delta_{nk}$ of **Equation (11)** is therefore evaluated for every segment by finding the shortest path on $G_{nk}(V, L_n, Y_{nk})$, and every segment attaining $\min_k \delta_{nk}$ is regarded as a possible assignment for record $n$. These segments yield feasible solutions of equal objective value, so the assignment among them is not uniquely identified by the observed tours.

**S4. Tour Extraction and Travel Cost Calculation**

The national origin–destination survey covered 1.18 million vehicles and reported the origin, destination, and purpose of each trip. In this study, passenger cars and chartered buses were targeted; taxis, cargo vehicles, and public buses were excluded. A vehicle was included if it made at least one leisure-related trip—sightseeing, tourism, sports, or recreational outings—and, for chartered buses, if the number of passengers remained unchanged at the destination, except on the first or last trip. For each such vehicle, consecutive leisure-related trips were linked into tours.

The travel link cost between each pair of municipal zones was calculated as follows (Armbrecht, 2014):

$$c_{na_1a_2} = \frac{1}{3}\alpha t_{a_1a_2} + \beta d_{a_1a_2} + f_{a_1a_2}, \tag{S1}$$

where $t_{a_1a_2}$, $d_{a_1a_2}$, and $f_{a_1a_2}$ denote the travel time, distance, and tolls, respectively, from municipal zone $a_1$ to $a_2$. The unit values $\alpha$= 48.8 [yen / min] and $\beta$= 16.3 [yen / km] are weighted averages of values published by Japan's Ministry of Land, Infrastructure, Transport and Tourism, weighted by the vehicle composition of the input data. The first and second terms of Equation (S1) represent the opportunity cost of travel time and fuel cost, respectively. Tolls were set to zero because they could not be identified in the survey.

## S5. Computational Performance

In the comparison examination, we sampled 20 datasets of 100 tours and set $|K| = 3$. Both methods were restricted to single-destination tours so that full enumeration remained feasible for the baseline. All computations in this study, including the main empirical analysis and this comparison, were performed using the Gurobi 10.0 optimizer on a computer equipped with fourteen 2.50 GHz processor cores and 128 GB memory, and the threshold $\theta$ in Algorithms S1 and S2 was set to $10^{-6}$.

The baseline algorithm explicitly formulates the complete problem by enumerating all candidate paths on the extended network and passes the resulting mixed-integer bilinear program directly to the solver. The proposed algorithm instead addresses the two computational difficulties separately: constraint generation adds only the violated shortest path constraints, and alternating optimization solves the restricted model by alternating between link-cost and segment-assignment updates.

The results are summarized in **Table S1**, which reports the average, minimum, and maximum values across the 20 datasets. Only a small subset of the alternatives was generated by the proposed algorithm, and no bilinear program was solved as a whole. The baseline's objective values are time-limit incumbents rather than proven optima.

**Table S1 Computational comparison of baseline and proposed algorithms**

| | Computation time [s] | | | Best objective value | | | Constraint-generation iterations | | |
|---|---|---|---|---|---|---|---|---|---|
| | avg. | min. | max. | avg. | min. | max. | avg. | min. | max. |
| Baseline algorithm | 3600.8 | 3600.6 | 3600.9 | 5091.8 | 1767.7 | 12778.7 | - | - | - |
| **Proposed algorithm** | **84.1** | **20.9** | **169.7** | **3947.7** | **0.0** | **26556.6** | **9.15** | **1** | **13** |

## S6. Comparison of Empirical Tour Optimality Gap Distributions with Gumbel Distributions

**Figure S2** shows the empirical cumulative distribution of the tour optimality gap for each segment, which is recovered rather than assumed and therefore takes whatever shape the observed tours imply. The distributions differ across segments: Segment 3C has a larger share of small gaps, whereas Segment 3A is more dispersed but reaches large gaps less often. For visual reference, the figure also plots Gumbel distributions with mode 0 and standard deviations of 1,000, 2,000, and 3,000. These curves are not fitted distributions or direct probabilistic counterparts of the recovered gaps. The comparison illustrates that the proposed method does not constrain the empirical gap distribution to a prespecified parametric form.

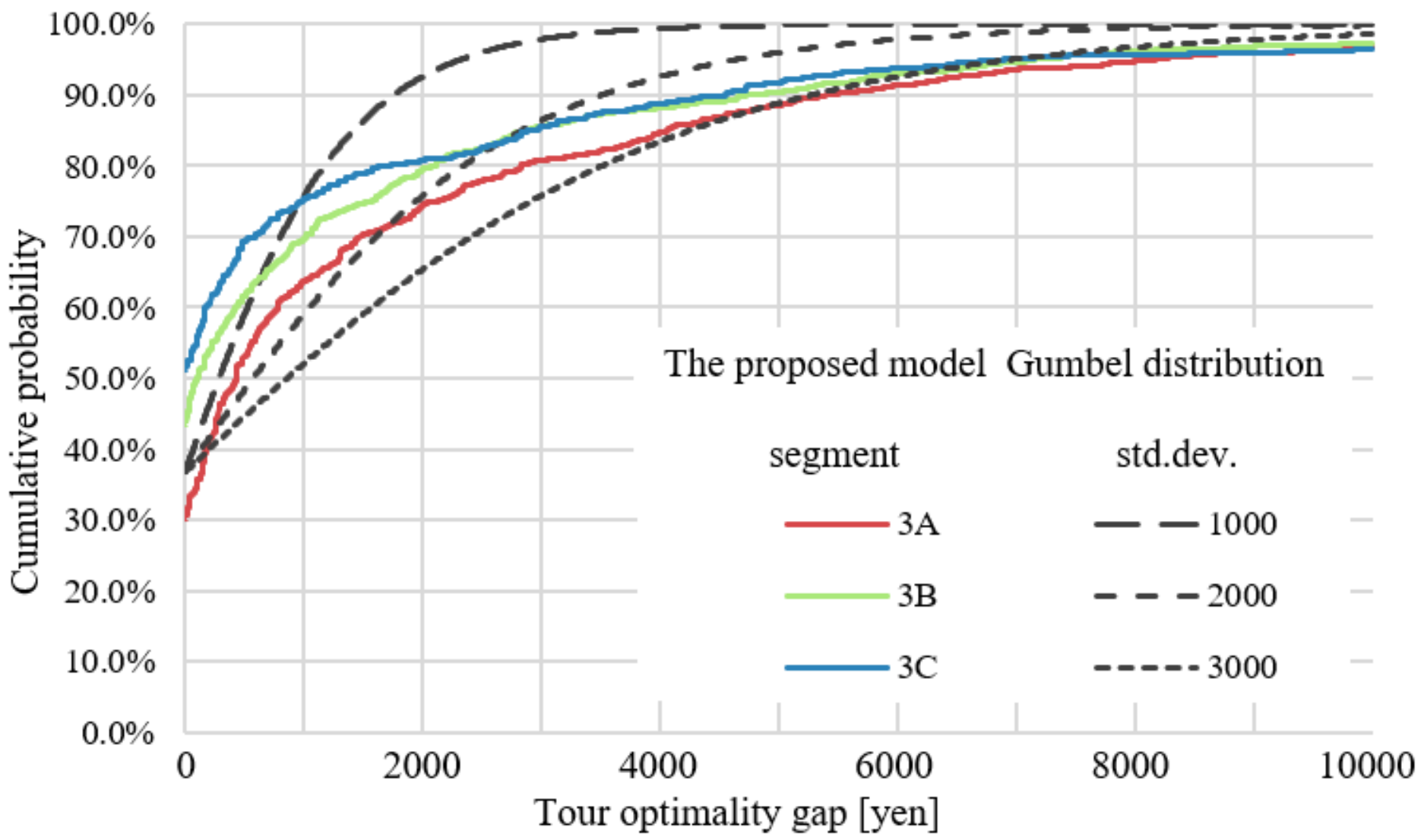

**Figure S2 Empirical cumulative distribution of the tour optimality gap by segment**

## References

Armbrecht J (2014) Use value of cultural experiences: A comparison of contingent valuation and travel cost. *Tourism Management* 42: 141–148.

Bishop CM (2006) *Pattern Recognition and Machine Learning*. New York: Springer.

Blei DM, Ng AY and Jordan MI (2003) Latent Dirichlet allocation. *Journal of Machine Learning Research* 3: 993–1022.

Cooper WW, Seiford LM and Zhu J (eds) (2011) *Handbook on Data Envelopment Analysis*. 2nd ed. New York: Springer.

Dempster AP, Laird NM and Rubin DB (1977) Maximum likelihood from incomplete data via the EM algorithm. *Journal of the Royal Statistical Society: Series B (Methodological)* 39(1): 1–22.

Mnih A and Salakhutdinov RR (2007) Probabilistic matrix factorization. In: Platt JC, Koller D, Singer Y and Roweis S (eds) *Advances in Neural Information Processing Systems 20*. Red Hook, NY: Curran Associates, pp.1257–1264.